\documentclass[11pt,twoside]{article}

\usepackage{amssymb}

\usepackage{amsmath}

\usepackage{amsfonts}

\usepackage{amssymb}

 \usepackage{bbm}
\allowdisplaybreaks[2]

\usepackage{latexsym}

\renewcommand{\skew}{\mathop{\rm skew}}
 \DeclareMathOperator{\sym}{sym}
\DeclareMathOperator{\tr}{tr} \DeclareMathOperator{\axl}{axl}
 \DeclareMathOperator{\dev}{dev}
\DeclareMathOperator{\Curl}{Curl\,} \DeclareMathOperator{\Div}{Div}
\DeclareMathOperator{\curl}{curl\,}
\DeclareMathOperator{\sL}{\mathfrak{sl}}
\DeclareMathOperator{\so}{\mathfrak{so}}

\newcommand{\Chi}{\raisebox{0.5ex}{\mbox{{\Large $\chi$}}}}

\newcommand{\yieldlimit}{\sigma_{\mathrm{y}}}

\newcommand{\C}{\mathbb{C}}

\newcommand{\BBR}{\mbox{$\mathbb{R}$}}

\newcommand{\SFH}{\mbox{$\mathsf{H}$}}

\newcommand{\SFQ}{\mbox{$\mathsf{Q}$}}
\newcommand{\SFV}{\mbox{$\mathsf{V}$}}

\newcommand{\SFZ}{\mbox{$\mathsf{Z}$}}
\newcommand{\bzero}{\mbox{$\bf 0$}}
\newcommand{\bvarepsilon}{\mbox{$\bf\varepsilon$}}
\newcommand{\btau}{\mbox{\bf$\tau$}}
\newcommand{\bm}{\mbox{$\mathfrak m$}}
\newcommand{\dsize}{\displaystyle}
\renewcommand{\div}{\mathop{\rm div}\nolimits}
\newcommand{\ba}{\mbox{\boldmath{$a$}}}
\numberwithin{equation}{section}
\newcommand{\parent}[3]{\left #1 {#3} \right #2}
\newcommand{\graffe}[1]{\parent \{ \}{#1}}

\newcommand{\norm}[1]{\left\Vert {#1} \right\Vert} 
\newcommand{\bfig}[2]{\begin{figure}\begin{center}\begin{picture}(341.8,#2)(
#1,0)}
\newcommand{\efig}[2]{\end{picture}\caption{#2.}\lbl{#1}\end{center}
\end{figure}}

\newcommand{\la}{\langle}
\newcommand{\ra}{\rangle}
\newcommand{\id}{{\boldsymbol{\mathbbm{1}}}}

\newcommand{\be}{\begin{equation}}

\newcommand{\ee}{\end{equation}}


\makeatletter
 \let\@fnsymbol\@arabic
 \makeatother
\begin{document}
\vskip-3truecm 
\title{
\vspace{-1.in} {\Large Existence results in dislocation based
rate-independent isotropic gradient plasticity with kinematical
hardening and plastic spin: The case with symmetric local backstress
}}

\author{{\large Fran\c{c}ois Ebobisse{\footnote{Corresponding author, Fran\c{c}ois Ebobisse,
Department of Mathematics and Applied Mathematics, University of
Cape Town, Rondebosch 7700, South Africa, e-mail:
francois.ebobissebille@uct.ac.za }}\quad and\quad Patrizio
Neff{\footnote {Patrizio Neff, Head of Chair of Nonlinear Analysis
and Modelling, Faculty of Mathematics, Universit\"at Duisburg-Essen,
Campus Essen, Universit\"atsstrasse 2, 45141 Essen, Germany, e-mail:
patrizio.neff@uni-due.de, http://www.uni-due.de/mathematik/ag\b{
}neff}}\quad and\quad Daya Reddy{\footnote{Daya Reddy, Department of
Mathematics and Applied Mathematics and Centre for Research in
Computational and Applied Mechanics, University of Cape Town, 7701
Rondebosch, South Africa, e-mail: daya.reddy@uct.ac.za,
http://www.mth.uct.ac.za/$\sim$bdr}}}\vspace{1mm}}



\date{\today}

 \maketitle

\begin{center}

\vspace{3mm}

\hspace{.05in}\parbox{5.5in} {{\small
\begin{abstract} In this paper we use convex analysis and variational inequality methods to establish an existence result for
a model of infinitesimal rate-independent gradient plasticity with
kinematic hardening and plastic spin, in which the local backstress
tensor remains symmetric. The model features a defect energy
contribution which is quadratic in the dislocation density tensor
$\Curl p$, giving rise to nonlocal non-symmetric kinematic
hardening. Use is made of
 a recently established Korn's type inequality for incompatible tensor fields. The solution space for the non-symmetric plastic distortion is
  naturally H(Curl) together with suitable tangential boundary conditions on the plastic distortion.
   Connections to other models are established as well.\end{abstract}}}
\end{center}
\noindent {\bf Key words:} plasticity, gradient plasticity,
dislocations, plastic spin, Korn's inequality, incompatible
distortions, rate-independent models, kinematical hardening,
backstress, variational inequality, defect energy.
\vskip.2truecm\noindent {\bf AMS 2010 subject classification:}
35B65, 35D10, 74C10, 74D10, 35J25, 75C05.

 {\small \tableofcontents}
\section{Introduction}\label{Intro} In the past twenty years, there have been
several experimental investigations for metallic and ceramic
materials which show that the elastic-plastic deformation of those
materials are size-dependent for sufficiently small scales. This
phenomenon cannot be predicted by conventional theories of
plasticity, which do not include any material length scales. Hence,
there clearly appeared a gap between micro-mechanical plasticity and
classical continuum plasticity. The purpose of the enhanced gradient
plasticity theories is to  formulate a constitutive framework on the
continuum level which is used to bridge that gap. In this paper we
will only discuss phenomenological models of isotropic
polycrystalline plasticity excluding the important case of single
crystal plasticity.

 There is now an abundant literature, and
research activities towards the development of models of gradient
plasticity which capture better the observed size-dependency
mentioned above. In the works of M\"uhlhaus and Aifantis
\cite{MUHAIF1991} and Gudmundson \cite{GUD2004}, the yield-stress is
set to depend also on some derivative of a scalar measure of the
accumulated plastic distortion. In the works of Gurtin and Anand
\cite{GURTAN2005}, Gudmundson \cite{GUD2004} and Neff et al.
\cite{NCA}, the yield-stress is not modified, but the free-energy is
augmented by a term involving the dislocation density. Also, it is
assumed in \cite{GURTAN2005} that the plastic flow is governed not
necessarily by the stress deviator (as in classical plasticity), but
more generally by microstress tensors that also satisfy a balance
law.

 In \cite{MUHAIF1991,GUD2004,GURTAN2005}, the plastic distortion
variable is assumed to be symmetric. However, as Gurtin and Anand
note \cite[p. 1626]{GURTAN2005}: "\ldots unless the plastic spin is
(explicitly) constrained to zero, constitutive dependencies on the
Burger tensor necessarily involve dependencies on the
(infinitesimal) plastic rotation". Note that even in classical
plasticity, the effect of plastic spin has been studied by several
authors like Dafalias \cite{Dafalias1984, Dafalias1985}, Mandel
\cite{Mandel1971, Mandel1973} and Kratochvil \cite{Kratochvil1971,
Kratochvil1973}, who were the first to suggest that a complete
macroscopic plasticity theory must include constitutive relations
involving also the plastic spin.

\vskip.2truecm On the other hand, though there are several theories
of gradient plasticity available in the literature, the results of
mathematical analysis for these problems are still rather scarce.
The first result of mathematical analysis on a model of gradient
plasticity was due to Djoko et al. \cite{DEMR1}. While the
developments by Gurtin-Anand \cite{GURTAN2005} were done for
viscoplastic bodies, the well-posedness of that model is considered
by Reddy at al. \cite{REM} for the rate-independent problem with
isotropic hardening and with both energetic and dissipative length
scales involved. Also, computational aspects of the model, based on
the work \cite{GURED2014}, are studied in \cite{BAREKLU2014} and are
devoted exclusively to single crystal plasticity. Let us mention
that, the purely energetic version of the Gurtin-Anand model, i.e.,
when the dissipative length scale $\ell=0$, is not yet treated but
will be done in this paper via the identification with the
irrotational version of Ebobisse-Neff \cite{EBONEFF} presented in
Paragraph \ref{irro-eboneff}. We also study in Section
\ref{Gurtin-kinlin}, the purely energetic case of the Gurtin-Anand
model with linear kinematic hardening. Another existence result for
the rate-independent problem of the Gurtin-Anand model was obtained
by
  Giaccomini and Lussardi \cite{giacoluss} within the
energetic-approach developed by Mielke \cite{Mielke2002, Mielke2003}
and it has also been proved that the model converges in a suitable
sense to a formulation of classical perfect plasticity proposed in
\cite{ddm} whenever the energetic
 and dissipative length scales go to zero.

Neff et al. proposed in \cite{NCA} a model of finite strain gradient
plasticity based on the multiplicative decomposition including
phenomenological Prager type symmetric linear kinematical hardening
and nonlocal kinematical hardening due to dislocations. The model is
from the outset non spin-free (the plastic distortion $p$ is not
symmetric) and its linearization leads to a thermodynamically
admissible model of infinitesimal plasticity involving only the Curl
of the non-symmetric plastic distortion $p$. The well-posedness of
the linear model is addressed as well, when formulated as a
variational inequality.

 In \cite{EBONEFF}, we have studied the
well-posedness of this model described within the framework of the
dual formulation for isotropic hardening by

 {\scriptsize\begin{table}[h!]
\begin{tabular}{|ll|}\hline &\qquad\\
{\em Additive split of distortion:} & $\nabla u=e+p$
\\
 {\em Additive split of strain:}& $\mbox{sym}\,\nabla u =\bvarepsilon=\bvarepsilon^e +\bvarepsilon^p$\,,\quad $\bvarepsilon^e=\sym\,e$\\
{\em Equilibrium:} & $\mbox{Div}\,\sigma +f=0$ with
$\sigma=\C.\bvarepsilon^e$\\&
 \\{\em Dissipation inequality:} &
 $\dsize\int_\Omega[\langle\sigma -\mu\,L_c^2\,\Curl\Curl p,\dot{p}\rangle -k_2\,\gamma\,\dot{\gamma}]
 \,dx\geq0$\\
 {\em Flow law in dual form:} & $\dot{p}\in\partial\Chi(\sigma-\mu \,L_c^2\,\mbox{Curl\,Curl\,}p),\qquad\dot{\gamma}=|\dot{p}|$\\
  & $\Chi$ is the indicator function of the set\\& of admissible
stresses\\
 \hline
\end{tabular}\caption{The model with plastic spin and isotropic hardening in \cite{EBONEFF}. }\end{table}}

\vskip.2truecm In this paper, we settle some of the open questions
raised in \cite{EBONEFF}. In that paper, we dealt with the isotropic
hardening case only and here we would like to extend our analysis to
the local linear Prager type kinematical hardening model which, in
the dual formulation of classical plasticity has the flow law
\begin{equation}\label{flowclassic}
\dot{\bvarepsilon}^p\in\partial\Chi(\sigma-b)\quad\mbox{ and
}\quad\dot{b}=\mu\, k_1\,\dot{\bvarepsilon}^p\,,\end{equation} in
which $b$ is the symmetric backstress tensor. Notice that
(\ref{flowclassic})$_2$ can be explicitly integrated (with proper
initial conditions $\bvarepsilon^p(0)=0$, $b(0)=0$) to yield
$b=\mu\, k_1\,\bvarepsilon^p$ and, substituted in
(\ref{flowclassic})$_1$, gives
$$\dot{\bvarepsilon}^p\in\partial\Chi(\sigma
-\mu\,k_1\,\bvarepsilon^p)\,.$$ Applying the same reasoning to our
nonlocal gradient model, we want to add a local backstress similar
to (\ref{flowclassic})$_2$. Hence, we define
$$
\dot{b}=\mu \,k_1\,\mbox{sym}\,\dot{p}\,,$$ in which we conserve the
symmetry of the local backstress tensor $b$. So, we get similarly,
$$
\dot{p}\in\partial\Chi(\sigma-b-\mu
L^2_c\,\mbox{Curl\,Curl}\,p)\quad\mbox{ and }\quad\dot{b}=\mu
\,k_1\,\mbox{sym}\,\dot{p}\,$$
 and integrating as above
 yields\begin{equation}\label{backstress-nonsym}
 \dot{p}\in\partial\Chi(\Sigma_E)\subset\BBR^+\frac{\dev\,\Sigma_E
}{|\dev\,\Sigma_E|}\quad\mbox{ with }\quad \Sigma_E=\sigma-\mu\,
k_1\,\mbox{sym}\,p-\mu\,
L_c^2\,\mbox{Curl\,Curl}\,p\,.\end{equation} Here, $\Sigma_E$ is the
elastic Eshelby tensor driving the plastic evolution. Note
immediately that it is only the nonlocal backstress contribution
$\mu\,L^2_c\,\Curl\Curl p$ which is responsible for the appearance
of plastic spin or not. In order to substantiate this claim, set
$L_c=0$ in (\ref{backstress-nonsym})$_2$ and consider
$$\dot{p}=\lambda\,\frac{\dev\,\Sigma_E}{|\dev\,\Sigma_E|}\in\mbox{Sym}(3)\,.$$ Assuming $p(0)\in\mbox{Sym}(3)$, we get that
$p(t)\in\mbox{Sym}(3)$ for every $t$ and hence we may replace $p$
with $\bvarepsilon^p:=\sym\,p$.

 From a purely mathematical point of view,
we could also consider the case of a non-symmetric local backstress
tensor, i.e.,
$$ \dot{p}\in\partial\Chi(\sigma-\widehat{b}-\mu\,
L^2_c\,\mbox{Curl\,Curl}\,p)\,,\quad\mbox{ and
}\quad\dot{\widehat{b}}=\mu\, k_1\,\dot{p}\,,$$
 which integrates to
 \begin{equation}\label{flowcurl}\dot{p}\in\partial\Chi(\sigma-\mu\,
 k_1\,p-\mu\,
L_c^2\,\mbox{Curl\,Curl}\,p)\,.\end{equation} In this case, the
following mathematical analysis would be drastically simplified
since no new estimate of the Korn type on incompatible tensor fields
(see \cite{NPW2014, NPW2011-1, NPW2012-1, NPW2012-2}) is needed. The
solution space is then trivially H(Curl). Notice that
(\ref{flowcurl}) also still reduces to a formulation of classical
plasticity in terms of a symmetric plastic strain tensor
$\bvarepsilon^p$ if the energetic length scale $L_c$ vanishes and
the initial plastic distortion $p(0)$ is chosen to be symmetric. To
see this, consider for $L_c=0$, the equation
$$\dot{p}=\lambda\,\frac{\dev\,\sigma
-\mu\,k_1\,p}{|\dev\,\sigma-\mu\,k_1\,p|}\,.$$ The format of the
equation, as far as classical solutions is concerned, is of the type
\begin{equation}\label{backstress-nonsym2}
\dot{p}=S(t)-\alpha\,p(t),\quad
p(0)=p_0\in\mbox{Sym}(3)\,,\end{equation} where
$S(t)\in\mbox{Sym}(3)$ and $\alpha\in\mathbb{R}$ can be assumed
given. Clearly, the system (\ref{backstress-nonsym2}) has only
symmetric solutions $p(t)$.\\
 The total energy would be of the type (isotropic elastic
response for simplicity)
$$\mu\,|\mbox{sym}\,(\nabla u-p)|^2 +\frac\lambda2\,|\mbox{tr}\,(\nabla u
-p)|^2+\underbrace{\frac{\mu\,
k_1}2\,|p|^2+\frac{\mu\,L^2_c}2\,|\mbox{Curl}\,p|^2}_{\mbox{\scriptsize{
immediate $H(\Curl)$-control  of $p$}}}\,,$$ which is, however, not
invariant w.r.t. the transformations $$u\to
u+\overline{A}x+\overline{b}\,,\qquad p\to p +\overline{A}\,,$$ for
constant skew-symmetric $\overline{A}\in\so(3)$ and constant
translation  $\overline{b}\in\mathbb{R}^3$, which represent
superposed Euclidean motions on both the displacement and plastic
distortion. Therefore, the choice $\mbox{sym}\,p$ in the backstress
evolution is mandatory by Euclidean invariance: the linear kinematic
hardening must be based on a symmetric backstress tensor.

Before we present our analysis of the model with linear kinematical
hardening and plastic spin, we find it important to first present,
using the convex analytical setting, a summary of those few models
of infinitesimal gradient plasticity in the literature for which a
mathematical analysis is now available. Precisely, we present in
Section \ref{various} the model by M\"uhlhaus-Aifantis
\cite{MUHAIF1991} as analyzed in \cite{DEMR1,DEMR2}, the model by
Gurtin-Anand \cite{GURTAN2005} as studied in \cite{REM}, the models
with plastic spin analyzed in \cite{NCA,EBONEFF} and their
irrotational version in \cite{NSW2009} and we highlight some
interconnnections. \vskip.2truecm\noindent
Let us first fix some
notations and definitions which will make the paper more clear and
readable.
\section{Some notational agreements and
definitions}\label{Notations} Let $\Omega$ be a bounded domain
 in $\BBR^3$ with Lipschitz continuous boundary $\partial\Omega$, which is occupied by an elastoplastic
body in its undeformed configuration. Let $\Gamma$ be a smooth
subset of $\partial\Omega$ with non-vanishing $2$-dimensional
Hausdorff measure. A material point in $\Omega$ is denoted by $x$
and the time domain under consideration is the interval $[0,T]$.\\
 For every $a,\,b\in\BBR^3$, we let $\la a,b\ra_{\BBR^3}$ denote the scalar
 product on $\BBR^3$ with associated vector
norm $|a|^2_{\BBR^3} = \la a, a\ra_{\BBR^3}$. We denote by
$\BBR^{3\times 3}$ the set of real $3\times 3$ tensors. The standard
Euclidean scalar product on $\BBR^{3\times 3}$ is given by $\la
A,\,B\ra_{\BBR^{3\times 3}} = \mbox{tr}\,\bigl[AB^T\bigr]$, where
$B^T$ denotes the transpose tensor of $B$. Thus, the Frobenius
tensor norm is $|A|^2 = \la A,\,A\ra_{\BBR^{3\times 3}}$.
 In the following we omit the subscripts $\BBR^3$ and $\BBR^{3\times 3}$. The identity tensor on $\BBR^{3\times 3}$ will be denoted by
  $\id$, so that $\mbox{tr}(A) = \la A, \id\ra$. We let
$\mbox{Sym\,}(3):=\{X\in\BBR^{3\times 3}\,|\,\,X^T=X\}$ denote the
set of symmetric tensors, the Lie-Algebras
$\so(3):=\{X\in\BBR^{3\times 3}\,|\,\,X^T=-X\}$ of skew-symmetric
tensors and  $\sL(3):=\{X\in\BBR^{3\times
3}\,|\,\,\mbox{tr\,}(X)=0\}$ of traceless tensors. For every
$X\in\BBR^{3\times 3}$, we set $\sym(X)=\frac12\bigl(X+X^T\bigr)$,
$\skew\,(X)=\frac12\bigl(X-X^T\bigr)$ and
$\dev(X)=X-\frac13\mbox{tr}\,(X)\,\id\in\sL(3)\,$ for the symmetric
part, the skew-symmetric part and the deviatoric part of $X$,
respectively.

The body is assumed to undergo infinitesimal deformations. Its
behaviour is governed by a set of equations and constitutive
relations. Below is a list of variables and parameters involved in
various models of infinitesimal gradient plasticity presented in
this paper:\begin{itemize}
\item[$\bullet$] $u$\,  the displacement of the macroscopic material
points;
\item[$\bullet$] $p$\, the plastic distortion variable is a
non-symmetric second order tensor, incapable of sustaining
volumetric changes; that is, $p\in\sL(3)$;

\item[$\bullet$] $e=\nabla u -p$\,  the elastic distortion  is a
non-symmetric second order tensor;

\item[$\bullet$] $\bvarepsilon^p=\sym\, p$\, the symmetric plastic strain
tensor;

\item[$\bullet$] $\bvarepsilon^e=\sym\,(\nabla u -p)$\, the symmetric elastic
strain tensor;

\item[$\bullet$] $\sigma$\, the Cauchy stress tensor is a symmetric
second order tensor;

\item[$\bullet$] $\yieldlimit$\, the yield stress;

\item[$\bullet$] $f$\, the body force;

\item[$\bullet$] $\Curl\,p=-\Curl\,e$\, the dislocation density
tensor;

\item[$\bullet$] $\tau^p$\, the microstress tensor is a second order
deviatoric symmetric tensor;

\item[$\bullet$] $\bm^p=(m_{ijk})$\, the micro-polar stress tensor is a third
order tensor deviatoric symmetric in the first two indices $i$ and
$j$. That is, $m_{ijk}=m_{jik}$ and $m_{iik}=0$;
\item[$\bullet$] $\gamma$\, the accumulated plastic strain.
\end{itemize}
\vskip.2truecm\noindent For isotropic media, the fourth order
elasticity tensor $\C$ is given by
\begin{equation}
\C. X = 2\mu\,\dev\,\sym X+\kappa \,\tr(X) \id =2\mu\,\sym
X+\lambda\,\tr(X)\id\label{C}
\end{equation}
for any second-order tensor $X$, where $\mu$ and $\lambda$ are the
Lam{\'e} moduli satisfying
\[
\mu>0\quad\mbox{ and }\quad 3\lambda +2\mu>0\,,
\] and $\kappa>0$ is the bulk modulus.\\
These conditions suffice for pointwise ellipticity of the elasticity
tensor in the sense that there exists a constant $m_0 > 0$ such that
\begin{equation}
\la X,\C. X\ra \geq m_0 |\sym X|^2\,. \label{ellipticityC}
\end{equation}

\vskip.2truecm\noindent The space of square integrable functions is
$L^2(\Omega)$, while the Sobolev spaces used in this paper are:
\begin{eqnarray}\label{sobolev-spaces}
&&\nonumber \mbox{H}^1(\Omega)=\{u\in
L^2(\Omega)\,|\,\mbox{grad\,}u\in
L^2(\Omega)\}\,,\qquad\quad\mbox{ grad}\,=\nabla\,,\\
&&\nonumber\qquad\quad
\norm{u}^2_{H^1(\Omega)}=\norm{u}^2_{L^2(\Omega)}+\norm{\mbox{grad\,}u}^2_{L^2(\Omega)}\,,\qquad\forall u\in\mbox{H}^1(\Omega)\,,\\
&&\nonumber \mbox{H}(\mbox{curl};\Omega)=\{v\in
L^2(\Omega)\,|\,\mbox{curl\,}v\in
L^2(\Omega)\}\,,\qquad\mbox{curl\,}=\nabla\times\,,\\
&&\nonumber\qquad\quad\norm{v}^2_{\mbox{\scriptsize
H}(\mbox{curl};\Omega)}=\norm{v}^2_{L^2(\Omega)}+\norm{\mbox{curl\,}v}^2_{L^2(\Omega)}\,,\,\,\quad\forall
v\in\mbox{H}(\mbox{curl;\,}\Omega)\,.
\end{eqnarray}
For every $X\in C^1(\Omega,\,\BBR^{3\times 3})$ with rows
$X^1,\,X^2,\,X^3$, we use in this paper the definition of $\Curl\,X$
in \cite{NCA, SVEN}:
$$\Curl X =\left(\begin{array}{l}\mbox{curl\,}X_1\,\,-\,\,-\\
\mbox{curl\,}X_2\,\,-\,\,-\\
\mbox{curl\,}X_3\,\,-\,\,-\end{array}\right)\in\BBR^{3\times 3}\,,$$
for which $\Curl\,\nabla v=0$ for every $v\in C^2(\Omega,\,\BBR^3)$.
Notice that the definition of $\Curl\,X$ above is such that $(\Curl
X)^Ta=\mbox{curl\,}(X^Ta)$ for every $a\in\BBR^3$ and this clearly
corresponds to the transpose of the Curl of a tensor as defined in
\cite{GURTAN2005, GURTAN-BOOK}.\\

In Paragraph \ref{pure-energ}, we will need an explicit definition
of the linear operator
$\mathbb{L}:\BBR^{3\times3\times3}\to\BBR^{3\times 3}$ such that
\begin{equation}\label{curl-grad}\Curl X=\mathbb{L}.\nabla X\qquad\forall X\in
C^1(\Omega,\,\BBR^{3\times 3})\,.\end{equation} So, for
$$
\overline{A}=\left(\begin{array}{ccc}
0 &-a_3&a_2\\
a_3&0& -a_1\\
-a_2& a_1&0
\end{array}\right)\in \so(3)\,,
$$
we consider the operator $\axl:\so(3)\rightarrow\mathbb{R}^3$
through
$$
\axl(\overline{A}):=\left( a_1, a_2, a_3 \right)^T,\quad \quad
\overline{A}.\, v=(\axl \overline{A})\times v,\quad \quad \forall \,
v\in\mathbb{R}^3,$$ $$ (\axl
\overline{A})_k=-\frac{1}{2}\,\sum\limits_{i,j=1}^3
\epsilon_{ijk}\overline{A}_{ij} =\frac{1}{2}\,\sum\limits_{i,j=1}^3
\epsilon_{kij}\overline{A}_{ji}\,,
$$
 where $\epsilon_{ijk}$ is the totally antisymmetric third order permutation
 tensor.\\
Hence, for every $A\in\BBR^{3\times 3}$, \begin{eqnarray*}
(\axl\skew A)_k=\frac{1}{2}\,\sum\limits_{i,j=1}^3
\epsilon_{kij}\skew(A)_{ji}&=&\frac{1}{4}\,\sum\limits_{i,j=1}^3
\epsilon_{kij}A_{ji}-\frac{1}{4}\,\sum\limits_{i,j=1}^3
\epsilon_{kij}A_ji\\
&=&\frac{1}{2}\,\sum\limits_{i,j=1}^3
\epsilon_{kij}A_{ji}\,.\end{eqnarray*}
 Recalling that $(\curl
v)_k=\dsize\sum_{i,j=1}^3\epsilon_{kij}v_{j,i}$ for every $v\in
C^1(\Omega,\BBR^3)$, it follows that
$$(\axl\skew\nabla v)_k=\frac{1}{2}\,\sum\limits_{i,j=1}^3
\epsilon_{kij}v_{j,i}=\frac12(\curl v)_k\,.$$
 Therefore, we may rewrite
$$
\Curl X=\left(
                       \begin{array}{c}
                         2\,\axl \skew \nabla X_1 \ \_\_\_ \\
                        2\,\axl \skew \nabla X_2 \ \_\_\_\\
                         2\,\axl \skew \nabla X_3 \ \_\_\_\\
                       \end{array}
                     \right)=\mathbb{L}. \nabla X,
$$
where $\mathbb{L}:\mathbb{R}^{3\times 3\times
3}\rightarrow\mathbb{R}^{3\times 3}$ is given by
\begin{equation}\label{operator-L}
\mathbb{L}=(\widetilde{\mathbb{L}}_1,\widetilde{\mathbb{L}}_2,
\widetilde{\mathbb{L}}_3)^T,\qquad \mathbb{L}. \nabla X=\left(
                       \begin{array}{c}
                         \widetilde{\mathbb{L}}_1. \nabla X \ \_\_\_ \\
                         \widetilde{\mathbb{L}}_2. \nabla X \ \_\_\_\\
                        \widetilde{\mathbb{L}}_3. \nabla X \ \_\_\_\\
                       \end{array}
                     \right)
\end{equation}
with $\widetilde{\mathbb{L}}_i:\mathbb{R}^{3\times
3}\rightarrow\mathbb{R}^{3}$, $i=1,2,3$ defined by
\begin{equation}\label{operator-Li} \widetilde{\mathbb{L}}_i. \nabla X=2\,\axl \skew
\nabla X_i.
\end{equation}
Hence, we have found the explicit linear operator
$\mathbb{L}:\mathbb{R}^{3\times 3\times
3}\rightarrow\mathbb{R}^{3\times 3}$, so that the equality
\eqref{curl-grad} holds.\\
Notice, that for every $X,\,Y\in C^1(\Omega,\,\BBR^{3\times 3})$
\begin{eqnarray*}
\la\Curl X,\Curl Y\ra&=&\sum_{i=1}^3\la\curl X_i,\curl Y_i\ra
=4\sum_{i=1}^3\la\axl\skew\nabla X_i,\axl\skew\nabla
Y_i\ra\\
&=&2\sum_{i=1}^3\la\skew\nabla X_i,\skew\nabla Y_i\ra =
2\sum_{i=1}^3\la\skew\nabla X_i,\nabla Y_i\ra\,.\end{eqnarray*}

The following function spaces and norm will be used later.
\begin{eqnarray}\label{Curl-spaces}
&&\nonumber \mbox{H}(\mbox{Curl};\,\Omega,\,\BBR^{3\times 3})=\{X\in
L^2(\Omega,\,\BBR^{3\times 3})\,|\,\mbox{Curl\,}X\in
L^2(\Omega,\,\BBR^{3\times 3})\}\,,\\
&&\,\norm{X}^2_{\mbox{\scriptsize H}(\mbox{\scriptsize
Curl};\Omega)}=\norm{X}^2_{L^2(\Omega)}+\norm{\mbox{Curl\,}X}^2_{L^2(\Omega)}\,,\quad\forall
X\in\mbox{H}(\mbox{Curl;\,}\Omega,\,\BBR^{3\times 3})\,,\\
&&\nonumber
\mbox{H}(\mbox{Curl};\,\Omega,\,\mathbb{E}):=\{X:\Omega\to\mathbb{E}\,|\,X\in
\mbox{H}(\mbox{Curl};\,\Omega,\,\BBR^{3\times 3})\}\,,
\end{eqnarray}
for $\mathbb{E}:=\mbox{Sym}\,(3),\,\sL(3)$ or
$\mbox{Sym}\,(3)\cap\sL(3)$.\vskip.2truecm\noindent
 We also consider the space
$$\mbox{H}_0(\mbox{Curl};\,\Omega,\,\Gamma,\BBR^{3\times 3})$$ as the completion in
the norm in  (\ref{Curl-spaces}) of the space $\{q\in
C^\infty(\Omega,\,\Gamma,\,\BBR^{3\times
3})\,|\,\,q\times\vec{n}|_\Gamma=0\}\,.$ Therefore, this space
generalizes the Dirichlet boundary condition
$$q\times\vec{n}|_\Gamma=0\,$$
to be satisfied by the plastic distortion $p$ or the plastic strain
$\bvarepsilon^p:=\sym p$. The space
$\mbox{H}_0(\mbox{Curl};\,\Omega,\,\Gamma,\mathbb{E})$ is defined as
 in (\ref{Curl-spaces}). \vskip.2truecm\noindent
 The divergence operator Div on second order
tensor-valued functions is also defined row-wise as
$$\mbox{Div}\,X=\left(\begin{array}{l}\mbox{div\,}X_1\\
\mbox{div\,}X_2\\
\mbox{div\,}X_3\end{array}\right)\,.$$

\section{Some models of infinitesimal gradient
plasticity}\label{various}
\subsection{The model by
M\"uhlhaus-Aifantis \cite{MUHAIF1991}} In this model, the
yield-stress in the case of isotropic hardening, is set to depend
also on some derivative of a scalar measure of the accumulated
plastic distortion which plays the role of the isotropic hardening
variable. A summary of the model is presented in Table
\ref{table:demr}.
 {\scriptsize\begin{table}[h!]
\begin{tabular}{|ll|}\hline &\qquad\\
 {\em Additive split of strain:}& $\mbox{sym}\,\nabla u =\bvarepsilon=\bvarepsilon^e +\bvarepsilon^p$\,,\quad $\bvarepsilon^p\in\mbox{Sym}\,(3)$\\
{\em Equilibrium:} & $\mbox{Div}\,\sigma +f=0$ with
$\sigma=\C.\bvarepsilon^e$\\ {\em Free energy:} &
$\frac12\langle\C.\bvarepsilon^e,\bvarepsilon^e\rangle+\frac12\,\mu\,k_2|\gamma|^2+
\frac12\,\mu\,l^2|\nabla\gamma|^2$\\ {\em Yield condition:} &
$\phi(\sigma,g)=|\dev\sigma|+g-\yieldlimit\leq0$\\
  where & $g=-\mu\,k_2\,\gamma +\mu\,l^2\Delta\gamma$\\&
\\{\em Dissipation inequality:} &
 $\dsize\int_\Omega\bigl[\langle\sigma,\dot{\bvarepsilon}^p\rangle+g\,\dot{\gamma}\bigr]dx\geq0$\\&\\
{\em Dissipation function}:
 &$\mathcal{D}(q,\xi):=\left\{\begin{array}{ll}\yieldlimit|q| &\mbox{ if
 }|q|\leq\xi,\\
 +\infty &\mbox{ otherwise}\end{array}\right.$
\\&\\
{\em Flow law in primal form:} & $(\sigma,g)\in\partial
\mathcal{D}(\dot{\bvarepsilon}^p,\dot{\gamma})$.\\
 {\em Flow law in dual form:} & $ \dot{\bvarepsilon}^p=\lambda\,\dsize\frac{\dev\sigma}{|\dev\sigma|}$,\quad $\dot{\gamma}=\lambda=|\dot{\bvarepsilon}^p|$\\&\\
 {\em KKT conditions:} & $\lambda\geq0$,\quad
 $\phi(\sigma,g)\leq0$,\quad
 $\lambda\,\phi(\sigma,g)=0$ \\
 {\em Boundary condition on $\gamma$:} & $\gamma=0$ on
 $\partial\Omega$\\ {\em Function spaces for $\bvarepsilon^p$ and $\gamma$:} & $\bvarepsilon^p(t,\cdot)\in L^2(\Omega,\,\mathbb{R}^{3\times 3})$,\quad $\gamma(t,\cdot)\in H^1_0(\Omega)$
  \\
 \hline
\end{tabular}\caption{The model by M\"uhlhaus-Aifantis
\cite{MUHAIF1991} as formulated in \cite{DEMR1,
DEMR2}\,.}\label{table:demr}\end{table}} \\
Under suitable boundary and initial conditions on
$u,\,\bvarepsilon^p$ and $\gamma$, the well-posedness as well as
computational aspects of that model are studied in \cite{DEMR1,
DEMR2} with the flow law formulated in its primal form.
\subsection{The model by Gurtin and Anand \cite{GURTAN2005} as
studied in \cite{REM} with isotropic
hardening}
 This model is based on the assumption that the power expended
by each kinematical field be expressible in terms of a system of
forces consistent with its own balance. Therefore, the model is
characterized by two additional stress tensors: a second order
tensor $\btau^p$ power conjugate to the symmetric plastic strain
$\bvarepsilon^p$ and a third order tensor $\bm^p$ power conjugate to
the gradient of the plastic strain, which satisfy a microforce
balance. The latter as well as the equilibrium being derived by the
principle of virtual power. Since $\bvarepsilon^p$ is deviatoric
symmetric, it is not restrictive to assume that $\tau^p$ is
deviatoric symmetric and the third order tensor $\bm^p$ is
deviatoric symmetric in the first two indices. The model as
formulated in \cite{REM} is summarized in Table \ref{table:rem} with
the purely energetic version in Table \ref{table:rem-energ}.

 {\scriptsize\begin{table}[h!]
\begin{tabular}{|ll|}\hline &\qquad\\
 {\em Additive split of strain:}& $\mbox{sym}\,\nabla u =\bvarepsilon=\bvarepsilon^e +\bvarepsilon^p$\,,\quad$\bvarepsilon^p\in\mbox{Sym}\,(3)$\\
{\em Equilibrium:} & $\mbox{Div}\,\sigma +f=0$ with
$\sigma=\C.\bvarepsilon^e$\\
{\em Microforce balance:} &$\dev\sigma=\btau^p-\mbox{Div}\,\bm^p$,\\
{\em where} & $\tau^p$: microstress ($2^{\mbox{\scriptsize nd}}$ order)\\
&$\bm^p$: micropolar stress ($3^{\mbox{\scriptsize rd}}$ order)
\\&\\ {\em Free energy:} &
$\frac12\langle\C.\bvarepsilon^e,\bvarepsilon^e\rangle+\frac12\mu\,
L^2_c\,|\mbox{Curl}\,\bvarepsilon^p|^2+\frac12\mu\,k_2|\gamma|^2$\\&\\
{\em Yield condition:} &
$\phi(\btau^p,\bm^p,g):=\sqrt{|\btau^p|^2+\ell^{-2}|\bm^p_{\rm\scriptsize
diss}|^2}+g-\yieldlimit\leq0$\\&\\
 {\em where } & $\bm^p_{\rm\scriptsize diss}=\bm^p-\bm^p_{\rm\scriptsize
 energ}$\\&\\
 &$\mu
 L^2_c\,\langle\mbox{Curl}\,\bvarepsilon^p,\mbox{Curl}\,\dot{\bvarepsilon}^p\rangle=\langle\bm^p_{\rm\scriptsize
 energ},\nabla\dot{\bvarepsilon}^p\rangle$\\&
 \\{\em Dissipation inequality:} &
 $\dsize\int_\Omega\bigl[\langle\tau^p,\dot{\bvarepsilon}^p\rangle+\langle\bm^p_{\rm\scriptsize
 diss},\nabla\dot{\bvarepsilon}^p\rangle+g\dot{\gamma}\bigr]dx\geq0$,\quad $g=-\mu\,k_2\,\gamma$\\&\\
 {\em Dissipation function:} &$\mathcal{D}(q,\xi):=\left\{\begin{array}{ll}\yieldlimit\, d^p(q) &\mbox{ if
 }d^p(q)\leq\xi,\\
 +\infty &\mbox{ otherwise}\end{array}\right.$\\&\\where & $d^p(q):=\sqrt{|q|^2+\ell^2|\nabla q|^2}$
\\ &\\{\em Flow law in primal form:} &$(\btau^p,\bm^p_{\rm\scriptsize diss},g)\in\partial
 \mathcal{D}(\dot{\bvarepsilon}^p,\nabla\dot{\bvarepsilon}^p,\dot{\gamma})$\\&\\
{\em Flow law in dual form:}
&$\left.\begin{array}{ll}\dot{\bvarepsilon}^p=\lambda\,\dsize\frac{\btau^p}{\yieldlimit
- g}, &
\nabla\dot{\bvarepsilon}^p=\lambda\,\dsize\ell^{-2}\frac{\bm^p_{\rm\scriptsize
diss}}{\yieldlimit-g},\\
 \dot\gamma=\lambda=d^p(\dot{\bvarepsilon}^p)&\end{array}\right\}\quad (*)$\\&\\
{\em KKT conditions:} &$\lambda\geq0$,\quad
$\phi(\btau^p,\bm^p_{\rm\scriptsize diss},g)\leq0$,\quad
$\lambda\,\phi(\btau^p,\bm^p_{\rm\scriptsize diss},g)=0$\\&\\
 {\em Boundary conditions for $\bvarepsilon^p$:} & $\bvarepsilon^p=0$ on
 $\partial\Omega$\\
 {\em Function space for $\bvarepsilon^p$:} & $\bvarepsilon^p(t,\cdot)\in H^1_0(\Omega,\,\mbox{Sym}\,(3))$\\
 {\em Two length scales:} & dissipative $\ell$ and energetic $L_c$\\
 \hline
\end{tabular}\caption{The model by Gurtin and Anand \cite{GURTAN2005} as formulated in \cite{REM} with dissipative {\bf
and} energetic length scales.}\label{table:rem}\end{table}}

 The well-posedness of the model was studied by Reddy et al.
\cite{REM}. We would like to emphasize here that the starting point
of the modelling and analysis is the primal form. The corresponding
solution will satisfy the dual form in which it is understood that
there is an extra consistency condition generated which makes
equation $(*)_2$ in Table \ref{table:rem} possible. The plastic
strain variable $\bvarepsilon^p$ is assumed from the outset to be
symmetric. Note that the formulation in \cite{GURTAN2005} as well as
in \cite{REM} involves the full gradient $\nabla\bvarepsilon^p$ of
the plastic strain in the dissipation function, which is controlled
in $L^2$ leading then to find the plastic strain variable
$\bvarepsilon^p$ in the Sobolev space
$H^1(\Omega,\,\mbox{Sym}\,(3))$ together with the possibility to
completely prescribe $\bvarepsilon^p$ at the boundary.

\subsection{The model with plastic spin in \cite{EBONEFF} and in \cite{NCA}}
Unlike the model in \cite{GURTAN2005} with the microstresses and the
plastic distortion kept symme\-tric, a model involving the plastic
spin is studied in \cite{NCA} with phenomenological Prager type
kinematical hardening and in \cite{EBONEFF} with isotropic
hardening. A summary of the setting in \cite{EBONEFF} for the
so-called equal spin case is presented in Table \ref{table:eboneff}.

{\scriptsize\begin{table}[h!]
\begin{tabular}{|ll|}\hline &\qquad\\
 {\em Additive split of distortion:}& $\nabla u =e +p$,\quad $\bvarepsilon^e=\mbox{sym}\,e$,\quad $\bvarepsilon^p=\sym p$\\
{\em Equilibrium:} & $\mbox{Div}\,\sigma +f=0$ with
$\sigma=\C.\bvarepsilon^e$\\&\\ {\em Free energy:} &
$\frac12\langle\C.\bvarepsilon^e,\bvarepsilon^e\rangle+\frac12\mu
\,L^2_c\,|\mbox{Curl}\,p|^2+\frac12\mu\,k_2\,|\gamma|^2$\\&\\
{\em Yield condition:} &
$\phi(\Sigma_E,g):=|\dev\Sigma_E|+g-\yieldlimit\leq0$\\&\\
 {\em where } & $\Sigma_E:=\sigma+\Sigma^{\mbox{\scriptsize
lin}}_{\mbox{\scriptsize curl}}$,\,\, $\Sigma^{\mbox{\scriptsize
lin}}_{\mbox{\scriptsize curl}}=-\mu \,L^2_c\,\Curl\Curl p$\\&\\ &
$g=-\mu\,k_2\,\gamma$\\&
 \\{\em Dissipation inequality:} &
 $\dsize\int_\Omega[\langle\Sigma_E,\dot{p}\rangle +g\dot{\gamma}]
 \,dx\geq0$\\&\\
 {\em Dissipation function:} &$\mathcal{D}(q,\xi):=\left\{\begin{array}{ll}\yieldlimit\, |q| &\mbox{ if }|q|\leq\xi,\\
 \infty & \mbox{otherwise}\end{array}\right.$\\&\\
 {\em Flow law in primal form:} &
 $(\Sigma_E,g)\in\partial \mathcal{D}(\dot{p},\dot{\gamma})$\\&\\
{\em flow law in dual form:}
&$\dot{p}=\lambda\,\dsize\frac{\dev\Sigma_E}{|\dev\Sigma_E|},\quad
\dot{\gamma}=\lambda=|\dot{p}|$\\&\\ {\em KKT conditions:}
&$\lambda\geq0$,\quad $\phi(\Sigma_E,g)\leq0$,\quad
$\lambda\,\phi(\Sigma_E,g)=0$\\&\\
 {\em Boundary conditions for $p$:} & $p\times\vec{n}=0$ on
 $\Gamma$,\,\, $(\Curl p)\times\vec{n}=0$ on $\partial\Omega\setminus\Gamma$\\
 {\em Function space for $p$:} & $p(t,\cdot)\in \mbox{H}(\mbox{Curl};\,\Omega,\,\BBR^{3\times 3})$\\
 \hline
\end{tabular}\caption{The models with plastic spin in \cite{EBONEFF}
and in \cite{NCA}\,.}\label{table:eboneff}\end{table}} \noindent An
existence result for the weak formulation of this model is obtained
in \cite{EBONEFF}. The solution space for this model is quite
naturally $p\in H(\Curl)$ since the isotropic hardening provides an
$L^2$-control of the entire plastic distortion $p$ and the energetic
 defect energy adds automatically a control of $\Curl p\in L^2$.
\subsection{The irrotational version of \cite{EBONEFF}.}\label{irro-eboneff} In \cite{NSW2009}, the irrotational limit case has been
computationally implemented as one of the first efficient treatments
of gradient plasticity. In this model, the plastic distortion $p$
remains symmetric and can therefore be written as
$\bvarepsilon^p=\sym p$. A summary of the model is presented in
Table \ref{table:irrot}. The well-posedness of this limit case is
included in the analysis presented in \cite{EBONEFF}.
{\footnotesize\begin{table}[h!]
\begin{tabular}{|ll|}\hline \\
 {\em Additive split of distortion:}& $\nabla u =e +p$,\quad $\bvarepsilon^e:=\mbox{sym}\,e$,\quad $\bvarepsilon^p:=\sym p$\\
&\\
{\em Equilibrium:} & $\mbox{Div}\,\sigma +f=0$ with $\sigma=\C.\bvarepsilon^e$\\&\\
{\em Free energy:} & $\frac12\langle\C.\bvarepsilon^e,\bvarepsilon^e\rangle+\frac12\mu \,L^2_c\,|\mbox{Curl}\,\bvarepsilon^p|^2+\frac12\mu\,k_2\,|\gamma|^2$\\&\\
{\em Yield condition:} &  $\phi(\Sigma_E,g):=|\dev\sym\Sigma_E|+g-\yieldlimit\leq0$\\
 & \\
 {\em where } & $\Sigma_E:=\sigma+\Sigma^{\mbox{\scriptsize
lin}}_{\mbox{\scriptsize curl}}$,\,\, $\Sigma^{\mbox{\scriptsize
lin}}_{\mbox{\scriptsize curl}}=-\mu \,L^2_c\,\Curl\Curl \bvarepsilon^p$\\&\\
& $g=-\mu\,k_2\,\gamma$
 \\&\\{\em Dissipation inequality:} &
 $\dsize\int_\Omega[\langle\Sigma_E,\dot{\bvarepsilon}^p\rangle +g\dot{\gamma}] \,dx\geq0$\\&\\
 {\em Dissipation function:} &$\mathcal{D}(q,\xi):=\left\{\begin{array}{ll}\yieldlimit |q| &\mbox{ if }|q|\leq\xi,\\\\
 \infty & \mbox{otherwise}\end{array}\right.$\\&\\

 {\em Flow law in primal form:} &
 $(\Sigma_E,g)\in\partial \mathcal{D}(\dot{\bvarepsilon}^p,\dot{\gamma})$\\&\\
{\em flow law in dual form:}
&$\dot{\bvarepsilon}^p=\lambda\,\dsize\frac{\dev\sym\Sigma_E}{|\dev\sym\Sigma_E|},\quad \dot{\gamma}=\lambda=|\dot{\bvarepsilon}^p|$\\&\\
{\em KKT conditions:} &$\lambda\geq0$,\quad
$\phi(\Sigma_E,g)\leq0$,\quad
$\lambda\,\phi(\Sigma_E,g)=0$\\&\\
 {\em Boundary conditions for $\bvarepsilon^p$:} & $\bvarepsilon^p\times\vec{n}=0$ on
 $\Gamma$,\,\, $(\Curl \bvarepsilon^p)\times\vec{n}=0$ on $\partial\Omega\setminus\Gamma$\\&\\
 {\em Function space for $\bvarepsilon^p$:} & $\bvarepsilon^p(t,\cdot)\in \mbox{H}(\mbox{Curl};\,\Omega,\,\mbox{Sym}(3))$\\&\\
 \hline
\end{tabular}\caption{The irrotational
version of \cite{EBONEFF} with isotropic hardening. The plastic
distortion
 itself does not appear, only $\bvarepsilon^p=\sym\,p$ remains in the model.}\label{table:irrot}\end{table}}

As shown in the next paragraph, this model can also be obtained as a
particular case of Gurtin-Anand \cite{GURTAN2005} for $l=0$,
$L_c>0$. Since the dissipative length scale $l=0$, the solution
space is only $H(\Curl)$ with the attendant tangential boundary
conditions. Thus, the existence result in \cite{EBONEFF} provides
also the first existence result for the purely energetic
Gurtin-Anand model with local isotropic hardening.

\subsection{The Gurtin-Anand model: purely energetic
version}\label{pure-energ}

In this section, we would like to compare or find a connection
between the model by Gurtin-Anand and our irrotational version.
 To this aim, we consider the defect energy
\begin{equation}\label{curlgrad}
\frac{1}{2}\mu\, L_c^2\, |\Curl \varepsilon^p|^2=\frac{1}{2}\mu\,
L_c^2\, |\mathbb{L}. \nabla \varepsilon^p|^2\,,
\end{equation}
where the  linear operator $\mathbb{L}:\mathbb{R}^{3\times 3\times
3}\rightarrow\mathbb{R}^{3\times 3}$ explicitly defined in
(\ref{operator-L})-(\ref{operator-Li}) is such that
$$\Curl \varepsilon^p=\mathbb{L}. \nabla
\varepsilon^p\,.$$ We recall that
$$
\Curl \varepsilon^p=\left(
                       \begin{array}{c}
                         \curl \varepsilon^p_1 \_\_\_ \\
                         \curl \varepsilon^p_2 \_\_\_\\
                         \curl \varepsilon^p_3 \_\_\_\\
                       \end{array}
                     \right),\quad \mbox{ $\varepsilon^p_i,\quad i=1,2,3$ denote the rows of
$\varepsilon^p\in \mathbb{R}^{3\times 3}$.}
$$

On the one hand, considering  the variation $\delta
\varepsilon^p_i\in C^\infty_0(\overline{\Omega},\Gamma)$ of the left
hand side of \eqref{curlgrad} with respect to the plastic strain
variable we get
\begin{eqnarray}\label{d1}
&&\frac{\rm d}{\rm dt}\frac12\int_\Omega \mu\, L_c^2  |\Curl
(\varepsilon^p+t \delta  \varepsilon^p)|^2dx\Big|_{t=0}= \mu\, L_c^2
\int_\Omega\langle \Curl \varepsilon^p, \Curl \delta
\varepsilon^p\rangle \, dx\notag\\
&&\quad=\,\mu\, L_c^2 \sum\limits_{i=1}^3\dsize\int_\Omega\langle
\curl \varepsilon^p_i, \curl \delta \varepsilon^p_i\rangle\, dx=
4\,\mu\, L_c^2 \sum\limits_{i=1}^3\int_\Omega\langle \axl \skew
\nabla \varepsilon^p_i, \axl \skew \nabla \delta
\varepsilon^p_i\rangle\, dx\notag\\
&&\quad=\,2\,\mu\, L_c^2 \sum\limits_{i=1}^3\int_\Omega\langle \skew
\nabla \varepsilon^p_i, \skew \nabla \delta \varepsilon^p_i\rangle\,
dx = \sum\limits_{i=1}^3\int_\Omega\langle \underbrace{2\,\mu\,
L_c^2 \skew \nabla \varepsilon^p_i}_{\mathfrak{m}^p_i,\  \text{2nd
order tensor}},
 \nabla \delta \varepsilon^p_i\rangle\, dx\notag\\
&&\quad=\,- \sum\limits_{i=1}^3\int_\Omega\langle
\Div\mathfrak{m}^p_i, \delta \varepsilon^p_i\rangle\, dx +
\sum\limits_{i=1}^3\int_{\partial\Omega}\langle  \mathfrak{m}^p_i.
\vec{n}, \delta \varepsilon^p_i\rangle\, da\notag
\\
&&\quad=\,- \sum\limits_{i=1}^3\int_\Omega\langle
\Div\mathfrak{m}^p_i, \delta \varepsilon^p_i\rangle\, dx
+\sum\limits_{i=1}^3\int_{\partial\Omega}\langle  (\axl
\mathfrak{m}^p_i)\times\vec{n},   \delta \varepsilon^p_i\rangle\,
da\notag
\\
&&\quad=\,-\sum\limits_{i=1}^3\int_\Omega\langle
\Div\mathfrak{m}^p_i, \delta \varepsilon^p_i\rangle\, dx -
\sum\limits_{i=1}^3\int_{\partial\Omega}\langle  \axl
\mathfrak{m}^p_i,   (\delta \varepsilon^p_i)\times \vec{n}\rangle\,
da\notag
\\
&&\quad=\,-\int_\Omega\langle  \Div\mathfrak{m}^p,   \delta
\varepsilon^p\rangle\, dx -
\sum\limits_{i=1}^3\int_{\partial\Omega}\langle  \axl
\mathfrak{m}^p_i,   (\delta \varepsilon^p_i)\times\vec{n}\rangle\,
da,
\end{eqnarray}
where the third order tensor $\mathfrak{m}^p\in \mathbb{R}^{3\times
3\times 3}$ is defined by
\begin{align}
\mathfrak{m}^p=(\mathfrak{m}^p_1,\mathfrak{m}^p_2,
\mathfrak{m}^p_3)^T.
\end{align}
On the other hand, we obtain
\begin{align}\label{d2}
\frac{\rm d}{\rm dt}\frac12\int_\Omega &\mu\, L_c^2  |\Curl
(\varepsilon^p+t \delta  \varepsilon^p)|^2dx\Big|_{t=0}
= \mu\, L_c^2 \int_\Omega\langle \Curl \varepsilon^p, \Curl \delta \varepsilon^p\rangle \, dx\notag\\
&=\mu\, L_c^2 \int_\Omega\langle \Curl \Curl \varepsilon^p, \delta \varepsilon^p\rangle \, dx
+\mu\, L_c^2 \sum\limits_{i=1}^{3}\int_{\partial \Omega}\langle  \delta \varepsilon^p_i\times [\Curl \varepsilon^p]_i, n\rangle \, da\\
&=\mu\, L_c^2 \int_\Omega\langle \Curl \Curl \varepsilon^p, \delta
\varepsilon^p\rangle \, dx-\mu\, L_c^2
\sum\limits_{i=1}^{3}\int_{\partial \Omega}\langle  [\Curl
\varepsilon^p]_i, \delta \varepsilon^p_i\times n\rangle \, da.\notag
\end{align}

In view of \eqref{d1} and \eqref{d2}, we obtain
\begin{align}
\mu\, L_c^2 \int_\Omega\langle \Curl \Curl \varepsilon^p, \delta
\varepsilon^p\rangle \, dx=-\int_\Omega\langle  \Div\mathfrak{m}^p,
\delta \varepsilon^p\rangle\, dx,
\end{align}
for all $\delta \varepsilon^p\in \mbox{H}_0(\Curl;\Omega;\Gamma)$,
i.e., for $\delta \varepsilon^p_i\times \vec{n}=0$.\\

 Since we can assume that
 $\bvarepsilon^p$ is trace free symmetric, so is
$\delta\bvarepsilon^p$ and we may equivalently write
$$\mu\,L_c^2\int_\Omega\la\mbox{dev}\,\sym\Curl\Curl\bvarepsilon^p,\delta\bvarepsilon^p\ra\,dx=
-\int_\Omega\la\mbox{Div\,}\bm^p,\delta\bvarepsilon^p\ra\,dx\,.$$Thus,
we get that
$$-\mbox{Div}\underbrace{\bm^p}_{\in\,
\mathbb{R}^{3\times3\times3}}=\underbrace{\mu\,L_c^2\,\mbox{dev}\,\sym\Curl\Curl
\bvarepsilon^p}_{\in\,\BBR^{3\times3}}\,.$$ Set
\begin{equation}\label{identification}\tau^p:=\dev\sigma +\mbox{Div\,}\bm^p=\underbrace{\dev\sigma
-\mu\,L_c^2\,\mbox{dev}\,\sym\Curl\Curl
\bvarepsilon^p}_{=\,\dev\sym\Sigma_E}\,.\end{equation}Hence,
$\tau^p=\dev\sym\Sigma_E$. \vskip.2truecm\noindent Notice that the
second order tensor $\mbox{Div}\,m^p$ is trace free. In fact, from
the bracket $\la \bm^p,\nabla\bvarepsilon^p\ra$, it holds that (we
may assume) $m^p_{ijk}=m^p_{jik}$ (since $\bvarepsilon^p$ is
symmetric) and we nay also assume that
\begin{equation}\label{trace:m} m^p_{iik}=0\,.\end{equation}
  Now, it is clear that
tr$(\mbox{Div}\,\bm^p)=m_{iik,k}=0$ from (\ref{trace:m}) and hence
$\mbox{Div}\,\bm^p$ is symmetric and trace free.

\vskip.2truecm\noindent Find a summary of this model in Table
\ref{table:rem-energ}.
 {\scriptsize\begin{table}[h!]
\begin{tabular}{|ll|}\hline &\qquad\\
 {\em Additive split of strain:}& $\mbox{sym}\,\nabla u =\bvarepsilon=\bvarepsilon^e +\bvarepsilon^p$\,,\quad$\bvarepsilon^p\in\mbox{Sym}\,(3)$\\
{\em Equilibrium:} & $\mbox{Div}\,\sigma +f=0$ with
$\sigma=\C.\bvarepsilon^e$\\
{\em Microforce balance:} &$\dev\sigma=\btau^p-\mbox{Div}\,\bm^p$,\\
{\em where} & $\tau^p$: microstress ($2^{\mbox{\scriptsize nd}}$ order)\\
&$\bm^p$: micropolar stress ($3^{\mbox{\scriptsize rd}}$ order)\\ &
$\tr(\Div\bm^p)=0$\\ & $\tau^p=\dev\sigma
+\Div\bm^p\in\sL(3)\cap\mbox{Sym}\,(3)$
\\&\\ {\em Free energy:} &
$\frac12\langle\C.\bvarepsilon^e,\bvarepsilon^e\rangle+\frac12\mu\,
L^2_c\,|\mbox{Curl}\,\bvarepsilon^p|^2+\frac12\mu\,k_2|\gamma|^2$\\&\\
{\em Yield condition:} &
$\phi(\btau^p,g):=|\btau^p|+g-\yieldlimit\leq0$\\&\\
 {\em where } & $\bm^p=\bm^p_{\rm\scriptsize
 energ}$\\&\\
 &$\mu
 L^2_c\,\langle\mbox{Curl}\,\bvarepsilon^p,\mbox{Curl}\,\dot{\bvarepsilon}^p\rangle=\langle\bm^p,\nabla\dot{\bvarepsilon}^p\rangle$\\&
 \\{\em Dissipation inequality:} &
 $\dsize\int_\Omega\bigl[\langle\tau^p,\dot{\bvarepsilon}^p\rangle+g\dot{\gamma}\bigr]\,dx\geq0$,\quad $g=-\mu\,k_2\,\gamma$\\&\\
 {\em Dissipation function:} &$\mathcal{D}(q,\xi):=\left\{\begin{array}{ll}\yieldlimit\, |q| &\mbox{
 if }
  |q|\leq\xi,\\
 +\infty &\mbox{ otherwise}\end{array}\right.$
\\ &\\{\em Flow law in primal form:} &$(\btau^p,g)\in\partial
 \mathcal{D}(\dot{\bvarepsilon}^p,\dot{\gamma})$\\&\\
{\em Flow law in dual form:}
&$\dot{\bvarepsilon}^p=\lambda\,\dsize\frac{\btau^p}{\yieldlimit -
g}=\dsize\frac{\tau^p}{|\tau^p|},\qquad
 \dot\gamma=\lambda=|\dot{\bvarepsilon}^p|$\\&\\
{\em KKT conditions:} &$\lambda\geq0$,\quad
$\phi(\btau^p,g)\leq0$,\quad
$\lambda\,\phi(\btau^p,g)=0$\\&\\
 {\em Boundary conditions for $\bvarepsilon^p$:} & $\bvarepsilon^p\times\vec{n}=0$ on
 $\partial\Omega$\\
 {\em Function space for $\bvarepsilon^p$:} & $\bvarepsilon^p(t,\cdot)\in \mbox{H}(\Curl;\,\Omega,\,\mbox{Sym}\,(3))$\\
 {\em Length scale:} &  energetic $L_c$\\
 \hline
\end{tabular}\caption{The irrotational model by Gurtin and Anand \cite{GURTAN2005} with no dissipative length scale i.e.
$\bm^p_{\rm\scriptsize diss}=0$\,. Since $\tau^p$ has been
identified with $\dev\sym\Sigma_E$ in (\ref{identification}), the
model coincides with the irrotational version of \cite{EBONEFF}.}
\label{table:rem-energ}\end{table}}

\vskip.3truecm\noindent

Let us now repeat the formulation of the model with spin in more
details, in its dual and primal setting for the paper to be rather
self-contained.

\section{The model with linear kinematical hardening and plastic
spin}\label{linkin-spin}
\subsection{Strong formulation} {\bf The balance equation.} The
conventional macroscopic force balance leads to the equation of
equilibrium
\begin{equation}
\div \sigma + f = \bzero\,. \label{equil}
\end{equation}
\vskip.1truecm\noindent {\bf Constitutive equations.} The
constitutive equations are obtained from a free energy imbalance
together with a flow law that characterizes plastic behaviour. Since
the model under study involves plastic spin, we consider an additive
decomposition of the displacement gradient $\nabla u$ into elastic
and plastic components $e$ and $p$ as mentioned in the notational
section, so that
\begin{equation}
\nabla u = e + p \label{displ-grad}
\end{equation}
\vskip.2truecm
 We consider here a free
energy of the form
\begin{eqnarray}\label{free-eng-kin}
\Psi(\nabla u,p,\Curl p): &=&\underbrace{\Psi^{\mbox{\scriptsize
lin}}_e(e)}_{\mbox{\small elastic energy}}\,\,
+\,\,\,\underbrace{\Psi^{\mbox{\scriptsize lin}}_{\mbox{\scriptsize
curl}}(\Curl p)}_{\mbox{\small defect
energy (GND)}}\\
\nonumber &+&\underbrace{\Psi^{\mbox{\scriptsize
lin}}_{\mbox{\scriptsize
 kin}}(p)}_{\mbox{\small linear kinematical hardening energy}}\,,
 \end{eqnarray} where
 \begin{eqnarray*}&&\Psi^{\mbox{\scriptsize
lin}}_e(e):=\frac12\la
\bvarepsilon^e,\C.\bvarepsilon^e\ra,\quad\Psi^{\mbox{\scriptsize
lin}}_{\mbox{\scriptsize curl}}(\Curl p):=\frac12\mu\, L_c^2\,|\Curl
p|^2\,\mbox{ and }\,\\\\
&&\Psi^{\mbox{\scriptsize lin}}_{\mbox{\scriptsize
 kin}}(p):=\frac12\mu\,k_1\,|\mbox{dev\,sym}\,p|^2\,.\end{eqnarray*}

 $L_c$ is the
energetic length scale and $k_1$ is the dimensionless  hardening
modulus. The defect energy is related to geometrically necessary
dislocations (GNDs) and the Burger's vector.

The local free-energy imbalance states that
\begin{equation}
\dot{\Psi} - \la\sigma,\dot{e}\ra - \la\sigma,\dot{p}\ra  \leq 0\
. \label{2ndlaw}
\end{equation}
Now we expand the first term, substitute (\ref{free-eng-kin}) and
get
\begin{equation}\label{exp-2ndlaw}
\la\mathcal{C}\bvarepsilon^e-\sigma,\dot{\bvarepsilon}^e\ra-\la\sigma,\dot{p}\ra+\mu\,
L_c^2\,\la\mbox{Curl}\,p,\mbox{Curl}\,\dot{p}\ra+\mu\,
k_1\,\la\mbox{dev\,sym}\,p,\dot{p}\ra\leq0\, ,
\end{equation}
which, using arguments from thermodynamics gives the elasticity
relation
\begin{equation}
\sigma = {\C}.\bvarepsilon^e=2\mu\, \sym(\nabla u-p)+\lambda\,
\tr(\nabla-p)\id \label{elasticlaw}
\end{equation}
and the reduced dissipation inequality
\begin{equation}\label{reduced-diss}-\la\sigma,\dot{p}\ra+ \mu\,
L_c^2\,\la\Curl p,\Curl \dot{p}\ra+\mu\,
k_1\,\la\mbox{dev\,sym}\,p,\dot{p}\ra\leq0.\end{equation} Now we
integrate (\ref{reduced-diss}) over $\Omega$ and  get
\begin{eqnarray}\nonumber0&\geq&\int_{\Omega}\Bigl[-\la\sigma,\dot{p}\ra+
\mu\, L_c^2\,\la\Curl p,\Curl \dot{p}\ra+\mu\,
k_1\,\la\mbox{dev\,sym}\,p,\dot{p}\ra\Bigr]\,dx\\
\nonumber&=&\int_{\Omega}\Bigl[-\la\sigma,\dot{p}\ra+\mu\,
L_c^2\,\la\Curl\,\Curl p,\dot{p}\ra+\mu\, k_1\,\la\mbox{dev\,sym}\,p,\dot{p}\ra\\
&& \qquad+\sum_{i=1}^3\mbox{div}\Bigl(\mu\,
L_c^2\,\frac{d}{dt}p^i\times(\Curl
p)^i\Bigr)\Bigr]\,dx\,.\label{reduced-diss2}
\end{eqnarray}
Using the divergence theorem we obtain
\begin{eqnarray}\label{reduced-diss3}
\nonumber&&\int_{\Omega}\left[\la-\sigma+\mu\, L_c^2 \,\Curl\,\Curl
p,\dot{p}\ra+\mu\, k_1\,\la\mbox{dev\,sym}\,p,\dot{p}\ra\right]\,dx\\
 &&\qquad +\sum_{i=1}^3\int_{\partial\Omega}\mu\,
L_c^2\,\langle\dot{p}^i\times(\Curl p)^i,\vec{n}\rangle da\leq0\, .
\end{eqnarray}

In order to obtain a dissipation inequality in the spirit of
classical plasticity, we assume that the infinitesimal plastic
distortion $p$ satisfies the so-called {\it linearized insulation
condition} \begin{equation}\sum_{i=1}^3\int_{\partial\Omega}\mu\,
L_c^2\,\langle\frac{d}{dt}p^i\times(\Curl p)^i,\vec{n}\rangle
da=0\,.\label{lin-sul}\end{equation} This condition is satisfied if
we
 assume for instance that the boundary is a perfect conductor. This means that the tangential component of $p$ vanishes on
$\partial\Omega$.
 In the context of dislocation dynamics these conditions express the requirement that there is no flux of the Burgers vector across a hard boundary.
  Gurtin and Anand \cite{GURTAN2005} introduce the following different types of boundary conditions for the plastic distortion%
\begin{align}
    \partial_t p\times\vec{n}|_{\Gamma_{\rm hard}}&=0 \quad \text{"micro-hard" (perfect conductor)} \notag \\
    \partial_t p|_{\Gamma_{\rm hard}}&=0 \quad \text{"hard-slip"}\\
    \Curl p\times\vec{n}|_{\Gamma_{\rm hard}}&=0 \quad \text{"micro-free"}\, . \notag
\end{align}
We specify a sufficient condition for the micro-hard boundary
condition, namely \begin{equation}\label{bc-plastic}
       p\times\vec{n}|_{\Gamma_{\rm hard}}=0
\end{equation}
and assume for simplicity only $\Gamma_{\rm
hard}=\partial\Omega=\Gamma$. Note that this boundary condition
constrains the plastic slip in tangential direction only, which is
what we expect to happen at $\Gamma_{\rm hard}$.\\
 Under (\ref{lin-sul}), we then
obtain the dissipation inequality
\begin{equation}\label{dissipation-ineq-kin}
\int_{\Omega} \la\sigma +\Sigma^{\mbox{\scriptsize
lin}}_{\mbox{\scriptsize curl}}+\Sigma^{\mbox{\scriptsize
lin}}_{\mbox{\scriptsize kin}},\dot {p}\ra\,dx\geq0 \,
,\end{equation} where $$\Sigma^{\mbox{\scriptsize
lin}}_{\mbox{\scriptsize curl}}:=-\mu \,L_c^2\,\Curl\,\Curl
p\quad\mbox{ and }\quad \Sigma^{\mbox{\scriptsize
lin}}_{\mbox{\scriptsize kin}}:=-\mu\, k_1\,\mbox{dev\,sym}\,p\, .$$\\
{\bf The flow law.}
  We consider a yield function
defined for every $\dsize\Sigma_{\mbox{\scriptsize
E}}:=\sigma+\Sigma^{\mbox{\scriptsize lin}}_{\mbox{\scriptsize
curl}}+\Sigma^{\mbox{\scriptsize lin}}_{\mbox{\scriptsize kin}}$ by
\begin{equation}
\label{yield-funct-kinematical} \phi_0(\Sigma_{\mbox{\scriptsize
E}}):= |\dev\Sigma_{\mbox{\scriptsize E}}| - \yieldlimit
\end{equation}
Here $\yieldlimit$ is the yield stress of the material. So the set
of admissible (elastic) generalized stresses is
\begin{equation}\label{admiss-stress-kin}
\mathcal{K}_0:=\left\{\Sigma_{\mbox{\scriptsize
E}}=\sigma+\Sigma^{\mbox{\scriptsize lin}}_{\mbox{\scriptsize
curl}}+\Sigma^{\mbox{\scriptsize lin}}_{\mbox{\scriptsize
kin}}\,|\,\,\phi_0(\Sigma_{\mbox{\scriptsize E}})
   \leq0\right\}\,.\end{equation}
The maximum dissipation principle gives the normality
law\begin{equation}\label{normalcone} \dot{p}\in
N_{\mathcal{K}_0}(\Sigma_{\mbox{\scriptsize E}}))\end{equation}
where $\dsize N_{\mathcal{K}_0}(\Sigma_{\mbox{\scriptsize E}})$
denotes the normal cone to $\mathcal{K}_0$ at
$\dsize\Sigma_{\mbox{\scriptsize E}}$, which is the set of
generalised strain rates $\dot{p}$ that satisfy
\begin{equation}
\la\overline{\Sigma} - \Sigma_{\mbox{\scriptsize E}},\dot{p}\ra \leq
0\ \quad \mbox{for all}\ \overline{\Sigma} \in {\mathcal{K}_0}\ .
\label{normality2}
\end{equation} Notice that $N_{\mathcal{K}_0}=\partial\Chi_0$ where
$\Chi_0$ denotes the indicator function of the set $\mathcal{K}_0$
and $\partial\Chi_0$ denotes the subdifferential of the function
$\Chi_0$.\\ Whenever the yield surface $\partial\mathcal{K}_0$ is
smooth at $\dsize\Sigma^{\mbox{\scriptsize p}}$ then
$$\dot{p}\in
N_{\mathcal{K}_0}(\Sigma_{\mbox{\scriptsize
E}}))\quad\Rightarrow\quad\exists\lambda\mbox{ such that }
\dot{p}=\lambda\,\frac{\dev\Sigma_{\mbox{\scriptsize
E}}}{|\dev\Sigma_{\mbox{\scriptsize E}}|}$$ with the Karush-Kuhn
Tucker conditions: $\lambda\geq0$, $\phi(\Sigma_{\mbox{\scriptsize
E}})\leq0$ and $\lambda\,\phi(\Sigma_{\mbox{\scriptsize E}})=0$\,.

Using convex analysis (Legendre-transformation) we find that
\begin{equation}\label{primalflowlaw}\dot{p}\in
\partial\Chi_0(\Sigma_{\mbox{\scriptsize E}})\quad\Leftrightarrow\quad \Sigma_{\mbox{\scriptsize E}} \in
\partial \Chi^*_0 (\dot{p})\, ,
\end{equation} where $\Chi^*_0$ is the Fenchel-Legendre dual of the function $\Chi_0$ denoted in this context by $\mathcal{D}$,
 the one-homogeneous dissipation function for rate-independent processes. That is,
\begin{equation}\label{dissp-function-kin}
\mathcal{D}(q)\,=\,\sup \graffe{\la\sigma+\Sigma^{\mbox{\scriptsize
lin}}_{\mbox{\scriptsize curl}}+\Sigma^{\mbox{\scriptsize
lin}}_{\mbox{\scriptsize kin}},q\ra \,\,|\,\
\phi_0(\sigma+\Sigma^{\mbox{\scriptsize lin}}_{\mbox{\scriptsize
curl}}+\Sigma^{\mbox{\scriptsize lin}}_{\mbox{\scriptsize
kin}})\leq0}\,=\,
      \yieldlimit\,|q|\,.
\end{equation} We get from the definition of the subdifferential ($\Sigma_{\mbox{\scriptsize E}} \in
\partial \Chi^*_0 (\dot{p})$) that,
\begin{equation}
\mathcal{D} (q) \geq \mathcal{D}(\dot{p}) +
\la\Sigma_{\mbox{\scriptsize E}},q - \dot{p}\ra\quad \mbox{for any
}q. \label{dissinequality}
\end{equation}
That is,
\begin{equation} \mathcal{D}(q)\geq \mathcal{D}(\dot{p})+\la\sigma +\Sigma^{\mbox{\scriptsize lin}}_{\mbox{\scriptsize
curl}}+\Sigma^{\mbox{\scriptsize lin}}_{\mbox{\scriptsize
kin}},q-\dot{p}\ra\mbox{ \,for any
}q\,.\label{dissinequality2}\end{equation}\\ {\bf Strong formulation
of the model.} To summarize, we have obtained the following strong
formulation for the model of infinitesimal gradient plasticity with
kinematic hardening and plastic spin. The goal is to find:
\begin{itemize}\item[(i)] the displacement $u\in \SFH^1(0,T;
\SFH^1_0(\Omega,{\Gamma},\mathbb{R}^3))$,
\item[(ii)] the infinitesimal plastic distortion $p$ with $\sym p\in
\SFH^1(0,T;L^2(\Omega, \sL(3)))$,\, $\Curl p\in \SFH^1(0,T;
L^2(\Omega,\BBR^{3\times 3}))$ and $ \Curl\Curl p \in \SFH^1(0,T;
L^2(\Omega,\BBR^{3\times 3}))$\end{itemize}
 such that the content of Table \ref{table:kinhardspin} holds.
{\footnotesize\begin{table}[h!]
\begin{tabular}{|ll|}\hline &\qquad\\
 {\em Additive split of distortion:}& $\nabla u =e +p$,\quad $\bvarepsilon^e:=\sym e$, $\bvarepsilon^p:=\sym p$\\
{\em Equilibrium:} & $\mbox{Div}\,\sigma +f=0$ with
$\sigma=\C.\bvarepsilon^e$\\&\\ {\em Free energy:} &
$\frac12\,\langle\C.\bvarepsilon^e,\bvarepsilon^e\rangle+\frac12\,\mu\,
L^2_c\,|\mbox{Curl}\,p|^2+\frac12\,\mu\, k_1\,|\dev\sym p|^2$\\&\\
{\em Yield condition:} &
$\phi(\Sigma_E):=|\dev\Sigma_E|-\yieldlimit\leq0$\\&\\
 {\em where } & $\Sigma_E:=\sigma+\Sigma^{\mbox{\scriptsize
lin}}_{\mbox{\scriptsize curl}}+\Sigma^{\mbox{\scriptsize
lin}}_{\mbox{\scriptsize kin}}$,\,\, \\&\\
&$\Sigma^{\mbox{\scriptsize lin}}_{\mbox{\scriptsize curl}}=-\mu\,
L^2_c\,\Curl\Curl p$,\, $\Sigma^{\mbox{\scriptsize
lin}}_{\mbox{\scriptsize kin}}=-\mu\, k_1\,\dev\sym p$
 \\&\\{\em Dissipation inequality:} &
 $\dsize\int_\Omega\langle\Sigma_E,\dot{p}\rangle dx\geq0$\\
 {\em Dissipation function:} &$\mathcal{D}(q):=\yieldlimit |q|$\\&\\
 {\em Flow law in primal form:} &
 $\Sigma_E\in\partial \mathcal{D}(\dot{p})$\\&\\
{\em Flow law in dual form:}
&$\dot{p}=\lambda\,\dsize\frac{\dev\Sigma_E}{|\dev\Sigma_E|},\quad\qquad \lambda=|\dot{p}|$\\&\\
{\em KKT conditions:} &$\lambda\geq0$, \quad $\phi(\Sigma_E)\leq0$,
\quad $\lambda\,\phi(\Sigma_E)=0$\\&\\
 {\em Boundary conditions for $p$:} & $p\times\vec{n}=0$ on
 $\Gamma$,\,\, $(\Curl p)\times\vec{n}=0$ on $\partial\Omega\setminus\Gamma$\\
 {\em Function space for $p$:} & $p(t,\cdot)\in \mbox{H}(\mbox{Curl};\;\Omega,\,\BBR^{3\times 3})$\\
 \hline
\end{tabular}\caption{The model with linear kinematical hardening and plastic
spin.  The boundary condition on $p$ necessitates at least $p\in
\mbox{H}(\mbox{Curl};\,\Omega,\,\BBR^{3\times 3})$.
However, whether this is the case will only be proven at the end of the paper.}\label{table:kinhardspin}\end{table}}\\
\subsection{Weak formulation of the model} Assume that the strong formulation has a solution $(u,p,\gamma)$.  Let $v\in \SFH^1(\Omega,\mathbb{R}^3)$
with $v_{|\Gamma}=0$. Multiply the equilibrium equation with
$v-\dot{u}$ and integrate in space to get
\begin{equation}\label{weak-eq1}
\int_{\Omega}\la\sigma,\nabla v-\nabla\dot{u}\ra dx=\int_\Omega
f(v-\dot{u})dx\, .
\end{equation}
Using the symmetry of the stress tensor $\sigma$ and the elasticity
relation we get
\begin{equation}\label{weak-eq2}
\int_{\Omega}\la\C.\sym(\nabla u-p),\mbox{sym}(\nabla
v-\nabla\dot{u})\ra dx=\int_\Omega f(v-\dot{u})dx\, .
\end{equation}
Now, we take any $q\in C^\infty(\overline{\Omega},\sL(3))$ such that
$q\times \vec{n}=0$ on $\Gamma$ and we integrate
(\ref{dissinequality2}) over $\Omega$, integrate by parts the term
with Curl\,Curl using the boundary conditions
$$(q-\dot{p})\times\vec{n}=0\mbox{ on }\Gamma$$ and get
\begin{eqnarray}
\label{weak-eq-kine}\nonumber\int_\Omega
\mathcal{D}^{\mbox{\scriptsize kin}}_0(q)\,dx &\geq&\int_\Omega
\mathcal{D}^{\mbox{\scriptsize kin}}_0(\dot{p})\,dx +\int_\Omega
\la\sigma +\Sigma^{\mbox{\scriptsize lin}}_{\mbox{\scriptsize
curl}}+\Sigma^{\mbox{\scriptsize lin}}_{\mbox{\scriptsize
kin}},q-\dot{p}\ra \,dx\\
\nonumber &\geq &\int_\Omega \mathcal{D}^{\mbox{\scriptsize
kin}}_0(\dot{p})\,dx+\int_\Omega \la\C.\sym(\nabla
u-p),\mbox{sym}(q-\dot{p})\ra\, dx\\
&& -\,\int_\Omega\la\,\mu\, L_c^2\,\Curl\,\Curl p+\mu\,
k_1\,\mbox{dev\,sym\,}p,q-\dot{p}\ra \,dx\\
\nonumber &\geq&\int_\Omega \mathcal{D}^{\mbox{\scriptsize
kin}}_0(\dot{p})\,dx +\int_\Omega\la \C.\sym(\nabla
u-p),\mbox{sym}(q-\dot{p})\ra\, dx\\\nonumber &&\nonumber-\,\mu\,
L_c^2\,\int_\Omega\la\Curl p,\Curl(q-\dot{p})\ra\, dx -\mu\,
k_1\,\int_\Omega\la\mbox{sym}\,p,q-\dot{p}\ra \,dx\, .
\end{eqnarray}
Adding (\ref{weak-eq-kine}) to the weak formulation of the
equilibrium in (\ref{weak-eq2}), we get that

\begin{eqnarray}\label{weak-form-kine}
&&\nonumber\int_{\Omega}\left.\Bigl[\la\C.\sym(\nabla
u-p),\mbox{sym}(\nabla v-q)-\mbox{sym}(\nabla\dot{u}-\dot{p})\ra+\mu
L_c^2\la\Curl p,\Curl(q-\dot{p})\ra\right.\\
&&\left. \quad\quad+\mu
k_1\la\mbox{sym}\,p,\mbox{sym}\,q-\mbox{sym}\,\dot{p}\ra\right.\Bigr]\,dx
+\int_\Omega \mathcal{D}^{\mbox{\scriptsize kin}}_0(q)\,dx
-\int_\Omega \mathcal{D}^{\mbox{\scriptsize
kin}}_0(\dot{p})\,dx\nonumber\\ &&\qquad\geq\int_\Omega
f(v-\dot{u})\,dx\qquad\forall (v,q)\,.\end{eqnarray}

\subsection{Existence result for the new formulation} To prove the
existence result for the weak formulation (\ref{weak-form-kine}), we
follow the abstract machinery developed by Han and Reddy in
\cite{Han-ReddyBook} for mathematical problems in classical
plasticity and used for instance in Djoko et al. \cite{DEMR1}, Reddy
et al. \cite{REM}, Neff et al. \cite{NCA}, Ebobisse-Neff
\cite{EBONEFF} for models of gradient plasticity. To this aim,
(\ref{weak-form-kine}) is written as the variational inequality of
the second kind: find $w=(u,p)\in \SFH^1(0,T;Z)$ such that $w(0)=0$
and
\begin{equation}\label{wf}
\ba(\dot{w},z-w)+j_0(z)-j_0(\dot{w})\geq \langle
\ell,z-\dot{w}\rangle\mbox{ for every } z\in Z\mbox{ and for a.e.
}t\in[0,T]\,,\end{equation} where $Z$ is a suitable Hilbert space to
be constructed later,
\begin{eqnarray}
&&\ba(w,z):=\int_{\Omega}\Bigl[\la\C.(\sym(\nabla u-p)),\sym(\nabla
v-q)\ra+\mu\, L^2\la\Curl p,\Curl q\ra\\
\nonumber &&\hskip3truecm+\mu\, k_1\la\sym p,\sym q\ra\Bigr]dx\,,
\label{bilin-kin}
\\\nonumber\\
&& j_0(z):=\int_\Omega
\mathcal{D}^{\mbox{\scriptsize kin}}_0(q)\,dx\,,\label{functional-kin}\\
&&\langle \ell,z\rangle:=\int_\Omega
fv\,dx\,,\label{lin-form}\end{eqnarray} for $w=(u,p)$ and $z=(v,q)$
in
$\SFZ$.\\
The Hilbert space $Z$ is constructed in such a way that the
functionals $\ba$, $j_0$ and $\ell$ satisfy the assumptions in the
abstract result in \cite[Theorem 7.3]{Han-ReddyBook}. The key issue
here is the coercivity of the bilinear form $\ba$ on $Z$. From the
structure of the bilinear form $\ba$ and the functional $j_0$, a
natural attempt for the space of infinitesimal plastic distortions,
is to consider the closure $\SFH_{\mbox{\scriptsize sym
}}(\Curl,\,\Omega,\Gamma;\sL(3))$ of the linear subspace
$$\{q\in C^\infty(\overline{\Omega},\BBR^{3\times
3})\,\,|\,\,\tr{q}=0,\,\,q\times\vec{n}=0\mbox{ on }\Gamma\}$$ with
respect to the norm
\begin{equation}\label{newnorm-q}\norm{q}_{\mbox{\scriptsize
sym,\,curl}}^2:=\norm{\sym
q}^2_{L^2}+\norm{\mbox{Curl\,}q}^2_{L^2}\, .
\end{equation}

Motivated by the well-posedness question for our model \cite{NCA,
EBONEFF}, Neff et al. \cite{NPW2011-1, NPW2012-1, NPW2012-2,
NPW2014}, derived a new inequality extending Korn's inequality to
incompatible tensor fields, namely there exist a constant
$C(\Omega)>0$ such that
\begin{align}
\label{incompatible_korn}
\forall \, p\in \SFH(\Curl;\,\Omega,\,\BBR^{3\times 3})\, | \quad & p\times\vec{n}|_\Gamma=0:   \\
  & \underbrace{\|p\|_{L^2(\Omega)}}_{\text{plastic distortion}}\le C(\Omega)\,
     \Big( \underbrace{\|\sym p\|_{L^2(\Omega)}}_{\text{plastic strain}}+
     \underbrace{ \|\Curl p\|_{L^2(\Omega)}}_{\text{dislocation density}} \Big)\, .\notag
\end{align}
Here, $\Gamma\subset\partial\Omega$ with full two-dimensional
surface measure
 and the domain $\Omega$ needs to be {\bf sliceable}, i.e. cuttable into finitely many simply connected subdomains with Lip\-schitz boundaries.
  The
inequality \eqref{incompatible_korn} expresses the important fact
that controlling the plastic strain $\sym p$ and the dislocation
density $\Curl p$ in $L^2(\Omega)$ gives a control of the full
plastic distortion $p$ in $L^2(\Omega)$ provided the correct
boundary conditions are specified: namely the micro-hard boundary
condition. Since in the sequel we assume that $\mathrm{tr}(p)=0$
(plastic incompressibility) the quadratic terms in the thermodynamic
potential provide a control of the right hand side in
\eqref{incompatible_korn}. So, setting:
\begin{eqnarray}
\SFV&=&\mathsf{H}^1_0(\Omega,{\Gamma},\mathbb{R}^3)=\{v\in
\mathsf{H}^1(\Omega,\mathbb{R}^3)\,|\,
v_{|\Gamma}=0\},\label{space-u}\\
 \SFQ&=&\mbox{H}_0(\mbox{Curl};\,\Omega,\,\Gamma,\sL(3))=\mathsf{H}_{\mbox{\scriptsize sym }}(\Curl;\,\Omega,\Gamma;\sL(3)),\label{space-p}\\
   \SFZ&=&\SFV\times \SFQ\,,\label{product-space}
\end{eqnarray} equipped with the norms
\begin{eqnarray}
&&  \norm{v}_V:=\norm{\nabla v}_{L^2},\label{norm-V}\quad\qquad
\norm{q}^2_{Q}:=\norm{\mbox{sym}\,q}^2_{L^2}+\norm{\mbox{Curl\,}q}^2_{L^2},\label{norm-Q}\\
&&\norm{z}^2_{Z}:=\norm{v}^2_{V} +\norm{q}^2_{Q}\quad\mbox{ for }
z=(v,q)\in \SFZ\,. \label{norm-Z}
\end{eqnarray}

Let us show that the bilinear form $\ba$ is coercive on $\SFZ$.
 Let therefore $z=(v,q)\in \SFZ$.
\begin{eqnarray*}
\ba(z,z)& \geq&  m_0\norm{\sym\,\nabla v-\sym q}^2_2 +\mu\,
L_c^2\norm{\Curl q}_2^2 +\mu\, k_1\norm{\sym q}_2^2 \mbox{ (from (\ref{ellipticityC}))}\\
&=&m_0\left[\norm{\sym \,\nabla v}^2_2 +\norm{\sym q}^2_2 -2\la\sym
\,\nabla v, \sym p\ra\right]\\
&&\qquad+\,\mu L_c^2\,\norm{\Curl q}_2^2
+\mu\, k_1\norm{\sym q}_2^2\\
&\geq &m_0\left[\norm{\sym\,\nabla v}^2_2 +\norm{\sym q}^2_2
-\theta\norm{\sym \,\nabla
v}_2^2-\frac1\theta\norm{\sym q}_2^2\right]\\
&&\qquad +\mu \,L_c^2\norm{\Curl q}_2^2+\mu\, k_1\norm{\sym
q}_2^2\quad\mbox{ (using Young's inequality)}\\
&=& m_0(1-\theta)\norm{\sym \,\nabla
v}^2_2+\left[m_0\Bigl(1-\frac1\theta\Bigr)+\mu\,
k_1\right]\norm{\sym q}_2^2+\mu \,L_c^2\norm{\Curl
q}_2^2.\end{eqnarray*}

 So, choosing $\theta$ such that
$\displaystyle\frac{m_0}{m_0+\mu k_1}< \theta<1$ and using Korn's
first inequality, we find a  positive constant
$C(m_0,\mu,k_1,L_c,\Omega)>0$ such that
$$\ba(z,z)\geq C\left[\norm{v}_V^2+\norm{\mbox{sym}\,q}^2_2
+\norm{\Curl q}^2_2 \right]=C\norm{z}^2_{Z}\quad\forall z=(v,q)\in
\SFZ\,,$$ which proves the coercivity of our bilinear form and the
inequality (\ref{incompatible_korn}) shows the equivalence
$\SFQ=\SFH_{\mbox{\scriptsize sym}}(\Curl,\Omega,\Gamma;\sL(3))$.

\section{The Gurtin-Anand model with linear kinematical hardening: purely energetic
version}\label{Gurtin-kinlin} \vskip.1truecm\noindent {\bf
Constitutive equations.}
\begin{equation}
\nabla u = e + p \quad \Rightarrow\quad \mbox{sym}\,\nabla
u=\bvarepsilon^e+\bvarepsilon^p\,.\label{displ-grad}
\end{equation}
\vskip.2truecm\noindent
 We consider here a free
energy of the form
\begin{eqnarray}\label{free-eng-gurt-kin}
\Psi(\bvarepsilon^e,\bvarepsilon^p,\Curl \bvarepsilon^p):
&=&\underbrace{\Psi^{\mbox{\scriptsize
lin}}_e(\bvarepsilon^e)}_{\mbox{\small elastic energy}}\,\,
+\,\,\,\underbrace{\Psi^{\mbox{\scriptsize lin}}_{\mbox{\scriptsize
curl}}(\Curl \bvarepsilon^p)}_{\mbox{\small defect
energy (GND)}}\\
\nonumber &+&\underbrace{\Psi^{\mbox{\scriptsize
lin}}_{\mbox{\scriptsize
 kin}}(\bvarepsilon^p)}_{\mbox{\small linear kinematical hardening energy}}\,,
 \end{eqnarray} where
 \begin{eqnarray*}&&\Psi^{\mbox{\scriptsize
lin}}_e(\bvarepsilon^e):=\frac12\la
\bvarepsilon^e,\C.\bvarepsilon^e\ra,\quad\Psi^{\mbox{\scriptsize
lin}}_{\mbox{\scriptsize curl}}(\Curl \bvarepsilon^p):=\frac12\mu\,
L_c^2\,|\Curl
\bvarepsilon^p|^2\,\,\,\mbox{ and }\,\\\\
&&\Psi^{\mbox{\scriptsize lin}}_{\mbox{\scriptsize
 kin}}(\bvarepsilon^p):=\frac12\mu\,k_1\,|\dev\bvarepsilon^p|^2\,.\end{eqnarray*}
Following the development in section \ref{linkin-spin}, the free
energy imbalance taking into account the boundary condition of the
plastic strain variable
$$\bvarepsilon^p\times\vec{n}|_\Gamma=0$$ leads to the dissipation
inequality
\begin{equation}\label{diss-ineq-gurt-kin}
\int_\Omega\la\Sigma_E,\dot{\bvarepsilon}^p\ra\,
dx\geq0\,,\end{equation}

where
$$\Sigma_E:=\sigma-\mu\,k_1\,\dev\bvarepsilon^p-\mu\,L_c^2\Curl\Curl\bvarepsilon^p\,.$$
 \vskip.1truecm\noindent {\bf The flow law}.
The set of generalized stresses is
$$\mathcal{K}:=\{\Sigma\in\mbox{Sym}\,(3)\,\,|\,\,\,\phi(\Sigma):=|\dev\Sigma|-\yieldlimit\leq0\}\,.$$ Hence,
 following \cite{EBONEFF}, we get the flow in dual form
\begin{equation}\label{normalcone} \dot{\bvarepsilon^p}\in
N_\mathcal{K}(\Sigma_E)\end{equation} where
$N_{\mathcal{K}}(\Sigma_E)$ denotes the normal cone to $\mathcal{K}$
at $\Sigma_E$, which in case of smoothness reads as
$$\dot{\bvarepsilon}^p=\lambda\,\frac{\dev\,\Sigma_E}{|\dev\,\Sigma_E|}$$
with $\lambda\geq0$, $\phi(\Sigma_E)\leq0$,
$\lambda\,\phi(\Sigma_E)=0$.\\
The flow law in its  primal formulation reads as
$$\Sigma_E\in\partial\mathcal{D}(\dot{\bvarepsilon}^p)\,.$$ That is,
\begin{equation}\label{flow-law-primal}
\mathcal{D}(q)\,\geq\,\mathcal{D}(\dot{\bvarepsilon}^p)+\la\Sigma_E,q-\dot{\bvarepsilon}^p\ra\quad\forall\,
q \in\mbox{Sym}\,(3)\,,\end{equation} where $\mathcal{D}$ is the
dissipation function defined as
$$\mathcal{D}(q):=\{\la\Sigma,q\ra\,\,|\,\,\Sigma\in\mathcal{K}\}=\yieldlimit\,|q|\quad\forall\,
q\in\mbox{Sym}\,(3)\,.$$
 {\bf Weak formulation of the model}. Now arguing as in Section \ref{linkin-spin} and also as in
 Ebobisse-Neff
\cite[Section 3]{EBONEFF}, we obtain a weak formulation of the model
in the form of a variational inequality
\begin{eqnarray}\label{VI-gurtin-linkin}
 &&\nonumber\int_{\Omega}\left.\Bigl[\la\C.(\sym\nabla
u-\bvarepsilon^p),(\sym\nabla
v-q)-(\sym\nabla\dot{u}-\dot{\bvarepsilon}^p)\ra+\mu\,
L_c^2\la\Curl \bvarepsilon^p,\Curl(q-\dot{\bvarepsilon}^p)\ra\right.\\
&&\left. \quad\quad+\mu\,
k_1\la\bvarepsilon^p,q-\dot{\bvarepsilon}^p\ra\right.\Bigr]\,dx
+\int_\Omega \mathcal{D}(q)\,dx -\int_\Omega
\mathcal{D}(\dot{\bvarepsilon}^p)\,dx\nonumber\\
&&\qquad\geq\int_\Omega f(v-\dot{u})\,dx\qquad\forall
(v,q)\,.\end{eqnarray}

The existence and uniqueness result for the variational inequality
is easily obtained in the spaces
\begin{eqnarray*}
&&
u\in\mbox{H}^1(0,T;\,\mbox{H}^1_0(\Omega,\,\Gamma,\,\mathbb{R}^3))\,,\\
&&\bvarepsilon^p\in\mbox{H}^1(0,T;\,\mbox{H}_0(\Curl,\Gamma;\mbox{Sym}(3)\cap\sL(3)))\,,\end{eqnarray*}
as in Section \ref{linkin-spin} through \cite[Theorem
7.3]{Han-ReddyBook}, following the coercivity on the space
$$Z:=\mbox{H}^1_0(\Omega,\,\Gamma,\,\mathbb{R}^3)\times
\mbox{H}_0(\Curl,\Gamma;\mbox{Sym}(3)\cap\sL(3))\,,$$ of the
bilinear form
$$\mathbf{a}(w,z):=\int_{\Omega}\Bigl[\la\C.(\sym\nabla
u-p),(\sym\nabla v-q)\ra+\mu\, L_c^2\la\Curl p,\Curl q\ra +\mu\,
k_1\la p,q\ra\Bigr]\,dx\,$$ for every $w=(u,p)$, $z=(v,q)$ in $Z$.\\
Note that since $\bvarepsilon^p$ is already trace-free and
symmetric, the coercivity of the bilinear form $\mathbf{a}$ does not
the new Korn's type inequality in \cite{NPW2011-1, NPW2012-1,
NPW2012-2, NPW2014}, and in (\ref{incompatible_korn}).

\section{The infinitesimal elastic micromorphic model}
The same total energy
\begin{eqnarray}\label{micromorphic1}
\nonumber\mathcal{E}(u,p)&=&\int_\Omega\left.\Bigl[\langle\C.\,\sym(\nabla
u-p),\sym(\nabla u-p)\rangle\right.\\
&&\qquad\quad \left.+\frac{\mu\,k_1}2|\dev\sym
p|^2+\frac{\mu\,L_c^2}2|\Curl p|^2-\langle
f,u\rangle\right.\Bigr]\,dx
\end{eqnarray} is the starting point for a two-field minimization
formulation
$$\mathcal{E}(u,p)\quad\rightarrow\quad\mbox{ min. w.r.t } (u,p)\,,$$ in the sense of a micromorphic model (\cite{NGMPR2014, NGMP2014}). \\
The relation of (\ref{micromorphic1}) to our plasticity formulation
 (\ref{backstress-nonsym}) is that in (\ref{micromorphic1}) the
micromorphic distortion $p$ is determined directly by a global
energy minimization instead of a plastic flow rule. The microbalance
equation is obtained as follows. The first variation of
(\ref{micromorphic1}) with respect to $p$ gives
\begin{eqnarray*}
&&\int_\Omega\bigl[\la\C.\sym(\nabla u-p),\sym\,\delta
p\ra+\mu\,k_1\la\dev\sym\,p,\delta
p\ra+\mu\,L^2_c\,\la\Curl\,\Curl\,p,\delta p\ra\bigr]dx\\
&&\quad=\,\int_\Omega\bigl[\la\sym\,\C.\sym(\nabla u-p), \delta
p\ra+\la\mu\,k_1,\dev\sym\,p,\delta
p\ra\\
&&\quad\qquad\quad+\,\la\mu\,L^2_c\,\Curl\,\Curl\,p,\delta
p\ra\bigr]dx=0\,.\end{eqnarray*} The "microbalance" is then
 of the form
 $$\mu\,L^2_c\,\Curl\,\Curl\,p=\overbrace{\sym\,\C.\sym(\nabla
 u-p)}^{{\mbox{\scriptsize Cauchy stress $\sigma$}}}-\mu\,k_1\,\dev\sym\,p\,.$$ The well-posedness of such
  a model has been shown in Neff et al. in \cite{NGMP2014}. Hence, we get in this model
  $$\mu\,L^2_c\,\Curl\,\Curl\,p=\sigma-\mu\,k_1\,\dev\sym\,p$$
  or $$
  0=\underbrace{\sigma-\mu\,k_1\,\dev\sym\,p-\mu\,L^2_c\,\Curl\,\Curl\,p}_{=\,\Sigma_E}\,$$
  with
  $$p\times\vec{n}|_\Gamma=0\quad\mbox{ and }\quad(\Curl
  p)\times\vec{n}|_{\partial\Omega\setminus\Gamma}=0\,,$$
  instead of a dual flow law $\dot{p}\in\partial\Chi(\Sigma_E)$ in plasticity.
\section{Conclusion}

 The development of the model with plastic spin is
straightforward and involves only the addition of a quadratic defect
energy. The boundary conditions on the plastic distortion are
consistent both from the physical and the mathematical point of
view. The departure from classical plasticity is minimal. Choosing a
symmetric local kinematical backstress evolution necessitates to use
a new Korn's type inequality for incompatible plastic distortions.
Contrary to the presented alternative models, in which the energetic
length scale $L_c$ has only a "passive" role in that necessary
estimates are already obtained from the dissipative length scale
$\ell$, in this model it is only the interplay between the energetic
length scale and the symmetric local backstress  which makes the
problem well-posed. By identifying the irrotational Gurtin-Anand
model with only energetic length scale as a special limit of our
model with spin, we have been able to provide an existence theorem
for that model for both the isotropic hardening case
(\cite{EBONEFF}), as well as the local backstress case (this paper
with the same considerations as in \cite{EBONEFF}). Moreover, our
derivation of the model avoids the introduction of certain
additional "micro force balances". Let us also mention that, the
introduction of the irrotationality constraint appears, in our
general framework with spin to be neither advantageous nor
necessary, but simplifies the analysis considerably.\vskip.2truecm

It remains to be seen if, in the dual formulation of our model with
spin one may consider isotropic hardening driven by a symmetrized
measure of accumulated plastic straining
\begin{equation}\label{iso-symp}
\dot{\gamma}=|\sym \dot{p}|\,.\end{equation} This would be
conceptionally pleasing since kinematical hardening could then
exclusively be related to the GND-distribution via the energetic
length scale $L_c$ (and assuming $k_1=0$) while the SSD-distribution
would be described by the accumulated plastic straining. We need to
remark that (\ref{iso-symp}) does not seem to satisfy the additional
assumption of {\it maximal dissipation}, making it unsuitable to be
considered in the primal formulation. However, it is
well-established that the equivalence of the primal and dual
formulation is not satisfied in general for gradient plasticity.

\vskip3truecm\noindent {\bf Acknowledgements:} The research of F.
Ebobisse has been supported by the National Research Foundation
(NRF) of South Africa through the Incentive Grant for Rated
Researchers and the International Centre for Theoretical Physics
(ICTP) through the Associateship Scheme. The first draft of this
work was written at Essen (Germany) in January 2015 while F.
Ebobisse was visiting the Faculty of Mathematics of the University
of Duisburg-Essen.

\end{document}